
%

\magnification=1200
\hsize=11.25cm
\vsize=18cm
\parskip 0pt
\parindent=12pt
\voffset=1cm
\hoffset=1cm



\catcode'32=9

\font\tenpc=cmcsc10
\font\eightpc=cmcsc8
\font\eightrm=cmr8
\font\eighti=cmmi8
\font\eightsy=cmsy8
\font\eightbf=cmbx8
\font\eighttt=cmtt8
\font\eightit=cmti8
\font\eightsl=cmsl8
\font\sixrm=cmr6
\font\sixi=cmmi6
\font\sixsy=cmsy6
\font\sixbf=cmbx6

\skewchar\eighti='177 \skewchar\sixi='177
\skewchar\eightsy='60 \skewchar\sixsy='60

\catcode`@=11

\def\tenpoint{%
  \textfont0=\tenrm \scriptfont0=\sevenrm \scriptscriptfont0=\fiverm
  \def\rm{\fam\z@\tenrm}%
  \textfont1=\teni \scriptfont1=\seveni \scriptscriptfont1=\fivei
  \def\oldstyle{\fam\@ne\teni}%
  \textfont2=\tensy \scriptfont2=\sevensy \scriptscriptfont2=\fivesy
  \textfont\itfam=\tenit
  \def\it{\fam\itfam\tenit}%
  \textfont\slfam=\tensl
  \def\sl{\fam\slfam\tensl}%
  \textfont\bffam=\tenbf \scriptfont\bffam=\sevenbf
  \scriptscriptfont\bffam=\fivebf
  \def\bf{\fam\bffam\tenbf}%
  \textfont\ttfam=\tentt
  \def\tt{\fam\ttfam\tentt}%
  \abovedisplayskip=12pt plus 3pt minus 9pt
  \abovedisplayshortskip=0pt plus 3pt
  \belowdisplayskip=12pt plus 3pt minus 9pt
  \belowdisplayshortskip=7pt plus 3pt minus 4pt
  \smallskipamount=3pt plus 1pt minus 1pt
  \medskipamount=6pt plus 2pt minus 2pt
  \bigskipamount=12pt plus 4pt minus 4pt
  \normalbaselineskip=12pt
  \setbox\strutbox=\hbox{\vrule height8.5pt depth3.5pt width0pt}%
  \let\bigf@ntpc=\tenrm \let\smallf@ntpc=\sevenrm
  \let\petcap=\tenpc
  \normalbaselines\rm}

\def\eightpoint{%
  \textfont0=\eightrm \scriptfont0=\sixrm \scriptscriptfont0=\fiverm
  \def\rm{\fam\z@\eightrm}%
  \textfont1=\eighti \scriptfont1=\sixi \scriptscriptfont1=\fivei
  \def\oldstyle{\fam\@ne\eighti}%
  \textfont2=\eightsy \scriptfont2=\sixsy \scriptscriptfont2=\fivesy
  \textfont\itfam=\eightit
  \def\it{\fam\itfam\eightit}%
  \textfont\slfam=\eightsl
  \def\sl{\fam\slfam\eightsl}%
  \textfont\bffam=\eightbf \scriptfont\bffam=\sixbf
  \scriptscriptfont\bffam=\fivebf
  \def\bf{\fam\bffam\eightbf}%
  \textfont\ttfam=\eighttt
  \def\tt{\fam\ttfam\eighttt}%
  \abovedisplayskip=9pt plus 2pt minus 6pt
  \abovedisplayshortskip=0pt plus 2pt
  \belowdisplayskip=9pt plus 2pt minus 6pt
  \belowdisplayshortskip=5pt plus 2pt minus 3pt
  \smallskipamount=2pt plus 1pt minus 1pt
  \medskipamount=4pt plus 2pt minus 1pt
  \bigskipamount=9pt plus 3pt minus 3pt
  \normalbaselineskip=9pt
  \setbox\strutbox=\hbox{\vrule height7pt depth2pt width0pt}%
  \let\bigf@ntpc=\eightrm \let\smallf@ntpc=\sixrm
  \let\petcap=\eightpc
  \normalbaselines\rm}
\catcode`@=12

\tenpoint



\catcode`\@=11
\def\pc#1#2|{{\bigf@ntpc #1\penalty \@MM\hskip\z@skip\smallf@ntpc%
	\uppercase{#2}}}
\catcode`\@=12

\def\pointir{\discretionary{.}{}{.\kern.35em---\kern.7em}\nobreak
   \hskip 0em plus .3em minus .4em }

\def\qed{\quad\raise -2pt\hbox{\vrule\vbox to 10pt{\hrule width 4pt
   \vfill\hrule}\vrule}}

\def\rem#1|{\par\medskip{{\it #1}\pointir}}

\def\vspace[#1]{\noalign{\vskip#1}}

\def\abstract#1{\vbox{\eightpoint\narrower\narrower 
\pc ABSTRACT|\pointir #1}}


\def\enslettre#1{\font\zzzz=msbm10 \hbox{\zzzz #1}}
\def\setN{{\enslettre N}}
\def\setZ{{\enslettre Z}}
\def\setQ{{\enslettre Q}}

\def\setF{{\enslettre F}}

\def\LHS{{\rm LHS}}
\def\RHS{{\rm RHS}}

\long\def\maskbegin#1\maskend{}


\def\section#1{\goodbreak\par\vskip .3cm\centerline{\bf #1}
   \par\nobreak\vskip 3pt }

\long\def\th#1|#2\endth{\par\medbreak
   {\petcap #1\pointir}{\it #2}\par\medbreak}

\def\article#1|#2|#3|#4|#5|#6|#7|
    {{\leftskip=7mm\noindent
     \hangindent=2mm\hangafter=1
     \llap{{\tt [#1]}\hskip.35em}{\petcap#2}\pointir
     #3, {\sl #4}, {\bf #5} ({\oldstyle #6}),
     pp.\nobreak\ #7.\par}}
\def\livre#1|#2|#3|#4|
    {{\leftskip=7mm\noindent
    \hangindent=2mm\hangafter=1
    \llap{{\tt [#1]}\hskip.35em}{\petcap#2}\pointir
    {\sl #3}, #4.\par}}
\def\divers#1|#2|#3|
    {{\leftskip=7mm\noindent
    \hangindent=2mm\hangafter=1
     \llap{{\tt [#1]}\hskip.35em}{\petcap#2}\pointir
     #3.\par}}


\def\frac#1#2{{#1\over #2}}

\def\bfu{{\bf u}}
\def\bfv{{\bf v}}
\def\bfa{{\bf a}}

\font\KFracFont=cmr12 at 20pt
\def\KK{\mathop{\lower 4pt\hbox{\KFracFont K}}\limits}
\def\Kadd#1{{\atop \;#1\;}}

\rightline{June 6, 2014}
\bigskip

\centerline{\bf Hankel continued fraction and its applications} 
\bigskip
\centerline{\sl Guo-Niu Han}
\footnote{}{\eightpoint
{\it Key words and phrases.} Hankel determinant, automatic proof, Hankel continued fraction, 
automatic sequence, Thue-Morse sequence, modulo $p$, Stieltjes algorithm, 
integer sequence, periodicity, regular paperfolding sequence, Stern sequence, irrationality exponent\par
{\it 2010 Mathematics Subject Classification.} 
05-04, 05A15, 11A55, 11B50, 11B85, 11C20, 11J82, 11T99, 11Y65, 15-04, 15A15, 15B33, 30B70.}
\bigskip
\bigskip
{\narrower\narrower
\eightpoint
\noindent
{\bf Abstract}.\quad
The Hankel determinants of a given power series $f$ can be evaluated by using 
the Jacobi continued fraction expansion of $f$. 
However the existence of the Jacobi continued fraction needs that all
Hankel determinants of $f$ are nonzero. We introduce
{\it Hankel continued fraction}, whose existene and unicity are guaranteed without any condition for the power series $f$. The Hankel determinants
can also be evaluated by using the Hankel continued fraction.

It is well known that the continued fraction expansion of a
quadratic irrational number is ultimately periodic.
We prove a similar result for power series. 
If a power series $f$ over a finite field satisfies a 
quadratic functional equation, then the Hankel continued fraction 
is ultimately periodic.
As an application, we derive the Hankel determinants of several automatic sequences,
in particular, the regular paperfolding sequence. 
Thus we provide an automatic proof of
a result obtained by Guo, Wu and Wen, which was conjectured by
Coons-Vrbik.

}

\section{1. Introduction} 

Let $\setF$ be a field and $x$ be a parameter. We identify a sequence ${\bf a}=(a_0, a_1, a_2, \ldots)$ over $\setF$ and its generating function 
$f=f(x)=a_0+a_1x+a_2x^2+\cdots\in\setF[[x]]$. Usually, $a_0=1$.
For each $n\geq 1$ and $k\geq 0$
the Hankel determinant of the series $f$ (or of the sequence $\bf a$)  is defined by 
$$
H_n^{(k)} (f) := \left|
\matrix{ a_k & a_{k+1} & \ldots & a_{k+n-1} \cr
a_{k+1} & a_{k+2} & \ldots & a_{k+n} \cr
\ \vdots \hfill & \ \vdots \hfill & \ddots &
\ \vdots \hfill \cr
a_{k+n-1} & a_{k+n} & \ldots & a_{k+2n-2} \cr} \right|\in\setF.
\leqno{(1.1)}
$$
Let $H_n(f):=H_n^{(0)}(f)$, for short;
the {\it sequence of the Hankel determinants} of $f$ is defined to be: 
$$H(f):=(H_0(f)=1, H_1(f), H_2(f), H_3(f), \ldots).$$

The Hankel determinants play an important role in the study of
the irrationality exponent of automatic number.
In 1998, Allouche, Peyri\`ere, Wen and Wen 
proved that all Hankel determinants of the Thue-Morse sequence
are nonzero [APWW]. 
Bugeaud [Bu11] was able to prove that the irrationality exponent of the Thue-Morse-Mahler number is equal to~2 by using APWW's result.
Using Bugeaud's method, several authors obtained the following results:
first,
Coons [Co13] who proved that the irrationality exponent of the sum of the
reciprocals of the Fermat numbers is~2; then, Guo, Wu and Wen who
showed that 
the irrationality exponents of the regular paperfolding numbers are exactly 2 [GWW].
However, the evaluations of the Hankel determinants
still rely on the method developed by Allouche, Peyri\`ere, Wen and Wen,
which consists of proving sixteen recurrence relations between determinants (see [APWW, Co13, GWW]).  
A combinatorial proof of the results by APWW and Coons about the Hankel determinants is derived by Bugeaud and the author [BH13].
In our previous paper [Ha13] short proofs of those results are presented by using
Jacobi continued fraction.
\medskip

The Hankel determinants of a given power series $f$ can be evaluated by using 
the Jacobi continued fraction expansion of $f$
(see, e.g., [Kr98, Kr05, Fl80, Wa48, Vi83, Ha13]). 
However the existence of the Jacobi continued fraction needs that all
Hankel determinants of $f$ are nonzero. In Section 2 we introduce
{\it Hankel continued fraction}, whose existene and unicity are guaranteed without any condition for the power series. The Hankel determinants
can also be evaluated by using the Hankel continued fraction (see Theorem 2.1). Let $p$ be a prime number and $\setF_p=\setZ/p\setZ$ be the finite field of size $p$. In Section 3 we prove the following result.

\proclaim Theorem 1.1.
Let $p$ be a prime number and $F(x)\in \setF_p[[x]]$ be a power series
satisfying the following quadratic functional equation
$$A(x)+B(x)F(x)+ C(x)F(x)^2 = 0, \leqno{(1.2)}$$
where $A(x), B(x), C(x)\in \setF_p[x] $ are three polynomials with one of the following conditions 
\smallskip
(i) $B(0)=1, \ C(0)=0,\ C(x)\not=0$;
\smallskip
(ii) $B(0)=1,\ C(x)=0$;
\smallskip
(iii) $B(0)=1,\ C(0)\not=0,\ A(0)=0$;
\smallskip
(iv) $B(x)=0,\  C(0)=1, \ A(x)= -(a_k x^k)^2 + O(x^{2k+1})$ for some $k\in \setN$ and $a_k\not=0$ when $p\not=2$.
\smallskip\noindent
Then, the Hankel continued fraction expansion of $F(x)$ exists and 
is ultimately periodic. Also, the Hankel determinant sequence $H(F)$ 
is ultimately periodic.

On the one hand, there is no similar result with traditional 
Jacobi continued fraction because of that its existence is not guaranteed, 
on the other hand, it is well known that the continued fraction expansion of a
quadratic irrational number is ultimately periodic.
Notice that the Hankel continued fraction and the Hankel determinant sequence in Theorem 1.1 can be {\it entirely} calculated by Algorithm 3.3.
By using Theorem~1.1 we derive the Hankel determinants of 
several automatic sequences.

%
\medskip

\proclaim Theorem 1.2.
For each pair of positive integers $a,b$, let 
$$
G_{a,b}(x)= {1\over x^{2^a}} \sum_{n=0}^\infty {x^{2^{n+a}} \over 1-x^{2^{n+b}}}\in\setF_2[[x]].
\leqno{(1.3)}
$$
Then $H(G_{a,b})$ is ultimately periodic.


A list of Hankel determinants for the special cases of Theorem~1.2 obtained by Algorithm 3.3 
is given in Corollary 4.1. 
When $a=b=0$, we then reprove Coons's Theorem [Co13].
The cases, where $(a,b)=(2,1), (2,0), (1,1)$, are obtained in [Ha13] by using
the Jacobi continued fraction expansion.
The case, where $a=0$ and $b=2$  was conjectured by 
Coons and Vrbik [CV12] and recently proved by Guo, Wu and Wen [GWW]
by using APWW's method.
The sequence $G_{0,2}$ is usually called {\it regular paperfolding sequence} [WiRP, Al87].
\medskip

An ultimately periodic sequence is written in contract form by using the star sign. For instance, the sequence $\bfa=(1, (3,0)^*)$ 
represents $(1,3,0,3,0,3,0,\ldots)$, that is, $a_0=1$ and $a_{2k+1}=3,
a_{2k+2}=0$ for each positive integer $k$.
Recall that the {\it Rudin-Shapiro sequence} $(u_n)$ is defined by 
$$
\cases{
u_0=0, \cr
u_{2n}=u_n, \quad
u_{4n+1}=u_n, \quad
u_{4n+3}=1-u_{2n+1}. \quad (n\geq 0) \cr
}
\leqno{(1.4)}
$$

\proclaim Proposition 1.3.
Let $(u_n)$ be the Rudin-Shapiro sequence and
$$
\leqalignno{
	f_1(x)&= \sum_{n\geq 0} u_{n+1} x^n;  \quad
	f_2(x)= \sum_{n\geq 0} u_{n+2} x^n;  \quad
	f_3(x)= \sum_{n\geq 0} u_{n+3} x^n.\cr
\noalign{Then,}
	H(f_1)&\equiv (1, 0, 0, 1, 0, 0, 1, 0, 0, 1, 1, 1, 0, 0, 0, 0, 1, 1 )^* \pmod 2;  \cr 
	H(f_2)&\equiv (1, 0, 1, 1, 0, 1, 1, 0, 1, 1, 1, 1, 0, 0, 0, 1, 1, 1 )^* \pmod 2;  \cr 
	H(f_3)&\equiv (1, 1, 0, 1, 1, 1, 1, 1, 1, 0, 1, 1, 0, 0, 1, 0, 1, 0 )^* \pmod 2.  \cr 
}
$$

\medskip

Recall that Stern's sequence $(a_n)_{n=0,1,\ldots}$ is defined by
(see [BV13, St58])
$$
\cases{
a_0=0,\quad a_1=1 \cr
a_{2n}=a_n,\quad a_{2n+1}=a_n+a_{n+1}. \quad (n\geq 1) \cr
}
$$
The twisted version of Stern's sequence $(b_n)$ is defined by
 (see [BV13, Ba10, Al12]) 
$$
\cases{
b_0=0,\quad b_1=1,\cr
b_{2n}=-b_n, \quad   b_{2n+1}=-(b_n+b_{n+1}). \quad (n\geq 1)
}
$$
Let
$$
S(x)=\sum_{n=0}^\infty a_{n+1}x^n
\quad\hbox{\rm and}\quad
B(x)=\sum_{n\geq 0} b_{n+1} x^n.
$$
be the generating function for Stern's sequence and twisted Stern's
sequence. 

\proclaim Proposition 1.4.
The Hankel determinants of the Stern's sequence and the twisted Stern's 
sequence verify the following relations
$$
H_n(S)/2^{n-2} \equiv
H_n(B)/2^{n-2} \equiv (0,0,1,1)^* \pmod 2.
$$

The proofs of Theorem 1.2 and Propositions 1.3-4 are given in Section~4. 
The results obtained in the paper about Hankel determinants
can be used for studying the irrationality exponent [BHWY].

\section{2. Hankel continued fractions} 
\def\sec{2}

Let ${\bf u}=(u_1, u_2, \ldots)$ and ${\bf v}=(v_0, v_1, v_2, \ldots)$
be two sequences.
Recall that the {\it Jacobi continued fraction} attached to $(\bfu, \bfv)$, 
or {\it $J$-fraction}, for short, is a continued
fraction of the form
$$
f(x)={v_0 \over 1 + u_1 x - 
	\displaystyle{v_1 x^2 \over 1+u_2x - 
		\displaystyle{v_2 x^2 \over {1 + u_3x - 
\displaystyle{v_3 x^2\over \ddots}}} }}, 
$$
The basic properties on $J$-fractions, we now recall, can be found in
[Kr98, Kr05, Fl80, Wa48, Vi83, Ha13].
The $J$-fraction of a given power series  $f$ exists  if and only if
all the Hankel determinants $H_n(f)$ are nonzero.
The first values of the coefficients $u_n$ and $v_n$ in the 
$J$-fraction expansion can be calculated by
the {\it Stieltjes Algorithm}.
Also, Hankel determinants can be calculated
from the $J$-fraction by means of the following  fundamental relation: 
$$
H_n(f)
= v_0^n v_1^{n-1} v_2^{n-2} \cdots v_{n-2}^2 v_{n-1}. 
$$

The Hankel determinants of a power series $f$ can be calculated by the 
above
fundamental relation 
if the $J$-fraction exists, which is equivalent to the fact that all Hankel determinants
of $f$
are nonzero.
In this section we define the so-called {\it Hankel continued fraction expansion} ({\it Hankel fraction} or {\it $H$-fraction}, for short)
whose existence and unicity are guaranteed without any condition for the power series.
The Hankel determinants can also be evaluated by using the Hankel continued fraction.

The relation between continued fractions and Hankel determinants
are widely studied. See [Kr05, Vi83, Fl80] for the $S$- and $J$-fractions;
[Bu10] and [Ci13] for $C$-fraction. 
The following table shows that the Hankel continued fraction
has some advantage over any other type of continued fractions.
$$
\def\tvi{\vrule height 12pt depth 5pt width 0pt}
\def\tv{\tvi\vrule}
\vbox{\offinterlineskip\halign{
\tv \ \hfil#\hfil\ \tv&&\strut\ \hfil#\hfil \cr   
\noalign{\hrule} 
		Fraction& Parameters  &\tv  & Fraction  &\tv& Fraction  &\tv    & Hankel det. &\tv   \cr
\noalign{\vskip -2.5mm}
		type &    &\tv  & existence  &\tv& unicity  &\tv    & formula &\tv   \cr
\noalign{\hrule} 
$S,J$-fraction&  $ \delta=1,2;  k_j=0$ &\tv  & No  &\tv& Yes  &\tv & Yes &\tv   \cr
$C$-fraction&  $ \delta=1, u_j(x)=0$ &\tv  & Yes  &\tv& Yes  &\tv & No &\tv   \cr
$H$-fraction & $\delta=2$   &\tv  & Yes  &\tv& Yes  &\tv & Yes &\tv   \cr
\noalign{\hrule}
}}
$$

{\it Definition \sec.1}. For each positive integer $\delta$, a {\it super continued fraction} associated with $\delta$, called {\it super $\delta$-fraction} for short, is defined to be a continued fraction of the following form
$$
F(x)
={v_0 x^{k_0}\over {
		1+u_1(x)x-\displaystyle{\strut v_1 x^{k_0+k_1+\delta}\over{
				1+	u_2(x)x-\displaystyle{\strut v_2 x^{k_1+k_2+\delta}\over{
						1+	u_3(x)x-\displaystyle{ \ddots
	}}}}}}}\leqno{(\sec.1)}
$$
where $v_j\not=0$ are contants, $k_j$ are nonnegative integers and $u_j(x)$ are polynomials of 
degree less than or equal to $k_{j-1}+\delta-2$. By convention, $0$ is of degree $-1$.
\medskip

When $\delta=1$ (resp. $\delta=2$) and all $k_j=0$, 
the super $\delta$-fraction (\sec.1)
is the traditional $S$-fraction (resp. $J$-fraction). 
A super $2$-fraction is called {\it Hankel continued fraction}.
When $\delta=1$ and $u_j(x)=0$, the super $1$-fraction is a special 
$C$-fraction (set $b_j=k_0+k_1+\cdots k_{j-1} + \lfloor j/2 \rfloor$ in [Ci13]).
Notice that every power series has a unique $C$-fraction expansion, but 
not all $C$-fractions have Hankel determinant formula,
and only those who are also super $1$-fractions have.

\proclaim Theorem \sec.1.
(i) Let $\delta$ be a positive integer.
Each super $\delta$-fraction defines
a power series, and conversely, for each power series $F(x)$,
the super $\delta$-fraction expansion of $F(x)$ exists
and is unique.
\smallskip
(ii) Let $F(x)$ be a power series such that its $H$-fraction 
is given by (\sec.1) with $\delta=2$. 
Then, all non-vanishing Hankel determinants of $F(x)$ are given by
$$
H_{s_j}(F(x))= (-1)^\epsilon v_0^{s_j} v_1^{s_j-s_1} v_2^{s_j-s_2} \cdots v_{j-1}^{s_j-s_{j-1}}, \leqno{(\sec.2)}
$$
where $\epsilon= \sum_{i=0}^{j-1} {k_i(k_i+1)/2}$
and
$s_j=k_0+k_1+\cdots + k_{j-1}+j$ for every $j\geq 0$.

The first part of Theorem \sec.1 is a consequence of Definition \sec.1 
and can be proved by using an algorithm. 
In fact, if $F(x)=v_0 x^{k_0} + O(x^{k_0+1})$ with $v_0\not=0$,
then, $F(x)/(v_0x^{k_0})=1+O(x)$. The polynomial $u_1(x)$
can be calculated by
$$
{v_0 x^{k_0} \over F(x)} = 1+ u_1(x)x - x^{k_0+\delta} F_1(x).
$$
We repeat the same operation for $F_1(x)$ and get $v_1, k_1, u_2(x)$,
etc.
The second part of Theorem \sec.1 follows from the next Lemma.

\proclaim Lemma \sec.2.
Let $k$ be a nonnegative integer and let $F(x), G(x)$ be two power series 
satisfying
$$
F(x)={x^k\over 1+u(x)x - x^{k+2} G(x)},\leqno{(\sec.3)}
$$
where $u(x)$ is a polynomial of degree less than or equal to $k$. Then,
$$
H_n(F)=(-1)^{k(k+1)/2} H_{n-k-1}(G). \leqno{(\sec.4)}
$$

{\it Proof}.
Let $F(x)=\sum_j f_jx^j$. 
We have $f_j=0$ for $j\leq k-1$ and $f_k=1$.
Let ${x^k/F(x)}=\sum_j b_j x^j$
and $G(x)=\sum_j g_j x^j$.
We have $g_j=-b_{j+k+2}$ for $j\geq 0$.
Let $b_k=f_j=0$ when $j<0$. 
We define four matrices by
$$
\leqalignno{
	{\bf F}_1&= (f_{i-j+k})_{0\leq i,j\leq n-1},\cr
	{\bf G}&= \hbox{\rm Diag} \bigl(
         (b_{i+j-k})_{0\leq i,j\leq k},\
	       (g_{i+j})_{0\leq i,j\leq n-k-1}\bigr),\cr
	{\bf F}&= (f_{i+j})_{0\leq i,j\leq n-1},\cr
	{\bf B}&= (b_{j-i})_{0\leq i,j\leq n-1},\cr
}
$$
and show that
$$
{\bf F}_1 \times {\bf G} = {\bf F} \times {\bf B}. \leqno{(\sec.5)}
$$
For example, when $k=3$, $n=7$, the four matrices and  (\sec.5) are
reproduced as follows.
$$
\leqalignno{
	&\left(\matrix{
1   &	.   &.   & .   &   . & . & .   \cr
f_4 &	1   &.   & .   &   . & . & .   \cr
f_5 &	f_4 &1   & .   &   . & . & .   \cr
f_6 &	f_5 &f_4 & 1   &   .   & . & .   \cr
f_7 &	f_6 &f_5 & f_4 &   1& . & .   \cr
f_8 &	f_7 &f_6 & f_5 &   f_4 & 1 & .   \cr
f_9 &	f_8 &f_7 & f_6 &   f_5 & f_4 & 1   \cr
}\right)\left(
\matrix{
.   &	.   &.   & 1   &   . & . & .   \cr
.   &	.   &1   & b_1   &   . & . & .   \cr
.   &	1   &b_1   &b_2   &   . & . & .   \cr
1   &	b_1 &b_2  & b_3   &   .   & . & .   \cr
. &	. &. & . &   g_0& g_1 & g_2   \cr
. &	. &. & . &   g_1 & g_2 & g_3   \cr
. &	. &. & . &   g_2 & g_3 & g_4   \cr
}
\right)\cr
&
\qquad
=
\left(\matrix{
.   &	.   &.   & 1   &   f_4 & f_5& f_6   \cr
.   &	.   &1   & f_4   &   f_5 & f_6 & f_7   \cr
.   &	1   &f_4   &f_5   &   f_6 & f_7 & f_8   \cr
1   &	f_4 &f_5  & f_6   &   f_7   & f_8 & f_9   \cr
f_4 &	f_5 &f_6 & f_7 &   f_8& f_9 & f_{10}   \cr
f_5 &	f_6 &f_7 & f_8 &   f_9 & f_{10} & f_{11}   \cr
		f_6 &	f_7 &f_8 & f_9 &   f_{10} & f_{11} & f_{15}   \cr
}
\right)
\left(\matrix{
1&b_1 & b_2 & b_3 &  b_4 & b_5 & b_6 	 \cr
.&1   & b_1 & b_2 &  b_3 & b_4 & b_5 	 \cr
.&.   &   1 & b_1 &  b_2 & b_3 & b_4 	 \cr
.&.   &  .  & 1   &  b_1 & b_2 & b_3 	 \cr
.&.   &   . & .   &  1   & b_1 & b_2 	 \cr
.&.   &   . & .   &  .   & 1   & b_1 	 \cr
.&.   &   . & .   &  .   & .   & 1   	 \cr
}\right)\cr
}
$$
Relations (\sec.5) are trivial for the entry $(i,j)$ 
when $0\leq j\leq k$ and when $j\geq k+1, i\leq k$.
For $i,j\geq k+1$. 
The two sides of (\sec.5) are
$$
\leqalignno{
	\LHS&=f_{i-1} g_{j-k-1} + f_{i-2} g_{j-k} + 
	\cdots + f_{i-n+k-1} g_{j+n-2k-1}\cr
	&=-(f_{i-1} b_{j+1} + f_{i-2} b_{j+2} + \cdots + f_{i-n+k-1} b_{j+n-k+1});\cr
	\RHS&=f_i b_j + f_{j+1} b_{j-1} \cdots + f_{i+n-1} b_{j-n+1}.\cr
}
$$
Since $F(x)\sum b_jx_j = x^k$, we have $\RHS-\LHS=0$. Moreover,
$\det{\bf F1}=1$, 
$\det{\bf G} = (-1)^{k(k+1)/2} H_{n-k-1}(G)$,
$\det{\bf F}=H_n(F)$, 
$\det{\bf B}=1$. This completes the proof of (\sec.4).\qed 

\medskip

{\it Example \sec.1}.
Let
$$f(x) =  {1- \sqrt{1-{4x^4\over 1+x}} \over 2x^4}\in\setQ[[x]].$$ 
Then
$$
f(x)={1\over {1+x-
\displaystyle{x^4\over 1 - 
\displaystyle{x^4 \over 1+x -
\displaystyle{x^4\over 1  - 
\displaystyle{x^4\over 1+x -
\displaystyle{x^4\over \ddots}}}  }}}}.
$$
In view of (2.1) we have $v_i=1$, $k_{2i}=0$, $k_{2i+1}=2$ for all $i$
and $(s_j)_{j=0,1,\ldots}=(0, 1, 4, 5, 8,9, 12,13,  \ldots)$
where $s_j$ is defined in Theorem \sec.1. By Theorem \sec.1 the Hankel determinant sequece is (see also [Ha13, Proposition 3.7])
$
H(f)=(1,1,0,0,-1,-1,0,0)^*.
$
\smallskip
{\it Example \sec.2}.
Let $g(x)$ be the generating function for the number of distinct partitions
$$
\leqalignno{
	g(x) &=  \prod_{n\geq1} (1+x^k) \in\setQ[[x]]\cr
	&= 1+x+x^2+2x^3+2x^4+3x^5+4x^6+5x^7+6x^8+8x^9+\cdots\cr
}
$$ 
Then
$$
g(x)={1\over {1-x-
\displaystyle{x^3\over 1 +x+ 
\displaystyle{x^5 \over 1-x+x^2-x^3 -
\displaystyle{x^5\over 1 +x+x^2+ 
\displaystyle{x^3\over 1-x +
\displaystyle{x^3\over \ddots}}}  }}}}.
$$
We have 
$$
\leqalignno{
	(k_j)_{j=0,1,\ldots}&= (0, 1, 2, 1, 0, 1, 0, 1, 0, 0, 0, 0, 0, 1, 1, 0, 0, 0, 2, \ldots),\cr
	(v_j)_{j=0,1,\ldots}&= (1, 1, -1, 1, -1, -1, -1, 1, -4, -1/4, 1/4, -8,  \ldots),\cr
	(s_j)_{j=0,1,\ldots}&= (0, 1, 3, 6, 8, 9, 11, 12, 14, 15, 16, 17, 18, 19,21, \ldots), \cr
	(H_j(g))_{j=0,1,\ldots}&= (1, 1, 0, -1, 0, 0, -1, 0, 1, 1, 0, -1, -1, 0, 1, -4,  \ldots).\cr
}
$$

\smallskip
{\it Example \sec.3}.
Let $
	h(x) =  (1-x)^{1/3} \in\setF_2[[x]].
$ 
Then
$$
h(x)={1\over {1+x+
\displaystyle{x^4\over 1 +x+ x^2+x^3+
\displaystyle{x^4 \over 1+x+ 
\displaystyle{x^8\over 1 +x+x^2+ x^3+
	\displaystyle{x^{16}\over\ddots}}  }}}}.
$$
We have 
$$
\leqalignno{
	(k_j)_{j=0,1,\ldots}&= (0, 2, 0, 6, 8, 22, 40 ,  \ldots),\cr
	(v_j)_{j=0,1,\ldots}&= (1, -1, -1, -1, -1, -1, -1 ,  \ldots),\cr
	(s_j)_{j=0,1,\ldots}&= (0, 1, 4, 5, 12, 21, 44, 85 ,  \ldots), \cr
	(H_j(h))_{j=0,1,\ldots}&= (1, 1, 0, 0, 1, 1, 0, 0, 0, 0, 0, 0, 1, 0, 0, 0, 0,  \ldots).\cr
}
$$


\section{3. The periodicity} 
\def\sec{3}

\medskip
In this section we prove Theorem 1.1. Let $\delta\in\setN^+$ and $\setF$ be a field. 

\proclaim Algorithm 3.1 [NextABC].
\smallskip\noindent
Prototype: $(A^*, B^*, C^*; k, A_k, D)=\hbox{\tt NextABC}(A,B,C; \delta)$
\smallskip\noindent
Input: $A(x), B(x), C(x)\in \setF[x]$ three polynomials such that 
$B(0)=1,$ $C(0)=0, C(x)\not=0, A(x)\not=0$;
\smallskip\noindent
Output: $A^*(x), B^*(x), C^*(x)\in\setF[x]$, $k\in \setN^+$, $A_k\not=0 \in\setF$,
$D(x)\in\setF[x]$ a polynomial of degree less than or equal to $k+\delta-1$ such that $D(0)=1$.
\smallskip\noindent
Step 1 [Define $k, A_k$]. Since $A(x)\not=0$, let $A(x)=A_kx^k + O(x^{k+1})$ with $A_k\not=0$.
\smallskip\noindent
Step 2. From (1.2) we have
$$
\leqalignno{
	F(x)&={-B+\sqrt{B^2-4AC}\over 2C}; &(3.1)\cr
	\noalign{\hbox{\rm and}}
	F(x)&={-A(x)\over B(x)+C(x)F(x)}. &(3.2)\cr
}
$$
Using (3.1) or (3.2) to get the first terms of $F(x)$, $F(x)/(-A_kx^k)$
and of $-A_kx^k/ F(x)$:
$$
\leqalignno{
	F(x)&=-A_kx^k + \cdots + O(x^{2k+\delta}); \cr
	{F(x)\over -A_kx^k}&=1 + \cdots + O(x^{k+\delta}); \cr
	{-A_kx^k\over F(x)}&=1 + \cdots + O(x^{k+\delta}). &(3.3)\cr
}
$$
\smallskip\noindent
Step 3 [Define $D$].  Define $D(x), G(x)$ by
$$
{-A_kx^k\over F(x)}= D(x)-x^{k+\delta} G(x) \leqno{(3.4)}
$$
where $D(x)$ is a polynomial of degree less than or equal to $k+\delta-1$ such that $D(0)=1$ and $G(x)$ is a power series. The value of $D(x)$ is obtained by (3.3).
\smallskip\noindent
Step 4 [Define $A^*, B^*, C^*$]. Let
$$
\leqalignno{
	A^*(x) &= \bigl(-D^2A/A_k+BD x^k-CA_kx^{2k}\bigr)/x^{2k+\delta}; \cr
	B^*(x) &= 2AD/(A_kx^{k}) - B ; &(3.5) \cr
	C^*(x) &= -Ax^{\delta}/A_k. \cr
}
$$
We prove that $A^*, B^*, C^*$ are polynomials in Lemma 3.2.

\smallskip

\proclaim Lemma 3.2.
Let $A(x), B(x), C(x)\in \setF[x]$ be three polynomials such that $B(0)=1, C(0)=0, C(x)\not=0, A(x)\not=0$ and  
$$(A^*, B^*, C^*; k, A_k, D)=\hbox{\tt NextABC}(A,B,C; \delta)$$
obtained by Algorithm 3.1.
If $F(x)$  is the power series defined by (1.2).
Then, 
$F(x)$ can be written as
$$
F(x) = {-A_k x^k\over D(x) - x^{k+\delta} G(x)} \leqno{(3.6)}
$$
where  $G(x)$ is a power series satisfying
$$A^*(x)+B^*(x)G(x)+ C^*(x)G(x)^2 = 0. \leqno{(3.7)}$$
Furthermore, $A^*(x), B^*(x), C^*(x)$ are three polynomials in $\setF[x]$ such that 
$B^*(0)=1, C^*(0)=0, C^*(x)\not=0$ and
$$
\deg(A^*)\leq d; \ 
\deg(B^*)\leq d+1; \ 
\deg(C^*)\leq d+\delta,   \leqno{(3.8)}
$$
where
$$
d=d(A,B,C)=\max(\deg(A)+\delta-2, \deg(B)-1, \deg(C)-\delta).
$$


{\it Proof}. From (1.2) and (3.6), we have
$$A(D - x^{k+\delta} G)^2+B{(-A_k x^k)( D - x^{k+\delta} G)}+ C {(-A_k x^k)}^2= 0.$$
Thus, $G(x)$ satisfies
$$\bar A(x)+\bar B(x)G(x)+ \bar C(x)G(x)^2 = 0 \leqno{(3.9)}$$
where 
$$
\leqalignno{
	\bar A &= AD^2-B A_k x^kD+CA_k^2 x^{2k}; \cr
	\bar B &= -2ADx^{k+\delta} + BA_k x^{2k+\delta}; \cr
	\bar C &= Ax^{2k+2\delta}. \cr
}
$$
Since $\bar C$ and $\bar B$ are divisible by $x^{2k+\delta}$,
so does $\bar A$.
Hence, (3.5) defines three polynomials $A^*, B^*, C^*$. Moreover,
$$
\leqalignno{
	\deg(A^*) &\leq \max(\deg(A)+\delta-2,  \deg(B)-1, \deg(C)-\delta); \cr
	\deg(B^*) &\leq \max(\deg(A)+\delta-1, \deg(B)); \cr
	\deg(C^*) &= \deg(A)+\delta. \cr
}
$$
Let $d_A=\deg(A),\ d_B=\deg(B)-1,\ d_C=\deg(C)-\delta$
and
$d_A^*=\deg(A^*)$,\ $d_B^*=\deg(B^*)-1,\ d_C^*=\deg(C^*)-\delta$.
The above inequalities become
$$
\leqalignno{
	d_A^* &\leq \max(d_A+\delta-2,  d_B, d_C); \cr
	d_B^* &\leq \max(d_A+\delta-2,  d_B); \cr
	d_C^* &= d_A. \cr
}
$$
So that $d_A^*, d_B^*, d_C^* \leq \max(d_A+\delta-2, d_B, d_C)$.
\qed

\proclaim Algorithm 3.3 [HFrac].
\smallskip\noindent
Prototype: $(a_k, d_k, D_k)_{k=0,1,\ldots} = \hbox{\tt HFrac}(A,B,C; p)$
\smallskip\noindent
Input: $p$ a prime number;
\smallskip\noindent
\qquad \quad $A(x), B(x), C(x)\in \setF_p[x]$ three polynomials such that $B(0)=1$, 
$C(0)=0$ and $C(x)\not= 0$;
\smallskip\noindent
Output: a finite or infinite sequence $(a_k, d_k, D_k)_{k=0,1,\ldots}$
\medskip\noindent
Step 1. $j:=0$,   $A^{(j)}:=A,\ B^{(j)}:=B,\ C^{(j)}:=C$.
\smallskip\noindent
Step 2. If $A^{(j)}= 0$, 
then return the finite sequence $(a_k, d_k, D_k)_{k=0,1,\ldots, j-1}$. The algorithm terminates.
\smallskip\noindent
Step 3. If $A^{(j)}\not= 0$, 
then let
$$(A^{(j+1)}, B^{(j+1)}, C^{(j+1)}; d_j, \alpha_j, D_j):=\hbox{\tt NextABC}(A^{(j)}, B^{(j)}, C^{(j)}; 2)).$$
Let $j:=j+1$.
\smallskip\noindent
Step 4. If there exits $0\leq i < j$ such that
$$(A^{(i)}, B^{(i)}, C^{(i)})=
(A^{(j)}, B^{(j)}, C^{(j)}),\leqno{(3.10)}$$
then return the infinite sequence 
$$((a_k, d_k, D_k)_{k=0,1,\ldots, i-1}, (a_k, d_k, D_k)_{k=i,i+1,\ldots, j-1}^* ).\leqno{(3.11)}$$
The algorithm terminates. Else, go to Step 2.
\smallskip\noindent

{\it Remarks}.
(i) In step 3 the conditions 
$$B^{(j)}(0)=1, C^{(j)}(0)=0, C^{(j)}(x)\not=0$$ are guaranteed by Lemma 3.2.
Algorithm~3.1 can be applied repeatedly. 
(ii) The loop Steps 2-4 will be broken at Step 2 or Step 4,
since the degrees of the polynomials $A^{(i)}, B^{(i)}, C^{(i)}$ are bounded, 
and the coefficients
are taken from $\setF_p$. The number of different
triplets $(A^{(i)}, B^{(i)}, C^{(i)})$ is finite. 
\medskip

\medskip
{\it Proof of Theorem 1.1}.
There are several cases to be considered. 
If $B(x)\not=0$, then we can always suppose that $B(x)=x^d+O(x^{d+1})$
for some $d\in \setN$.
\smallskip\noindent
(i) If $B(0)=1, C(0)=0, C(x)\not=0$, let 
$$(a_k, d_k, D_k)_{k=0,1,\ldots} = \hbox{\tt HFrac}(A,B,C; p).$$
By Lemma 3.2,
$$
F(x)=\displaystyle{ -a_0 x^{d_0} \over
	D_0(x)  + \displaystyle{a_1 x^{d_0+d_1+2}  \over
	D_1(x)  + \displaystyle{a_2 x^{d_1+d_2+2}  \over
	D_2(x)  + \displaystyle{a_3 x^{d_2+d_3+2}  \over \ddots
	}}}}
$$
and the above $H$-fraction is ultimately periodic (see Steps 2 and~4 in Algorithm 3.3). Note that if $A(x)=0$, then the output sequence 
$((a_k, d_k, D_k))$ $({k=0,1,2,\ldots})$ is the empty sequence. In this case $F(x)=0$.
\smallskip\noindent
(ii) If $B(0)=1, C(x)=0$, then $F(x)=-A(x)/B(x)$ is rational. 
\smallskip\noindent
(iii) If $B(0)=1, C(0)\not=0$ and $A(x)=0$, then $F(x)$ is rational.
If $B(0)=1, C(0)\not=0$ and $A(x)=A_k x^k + O(x^{k+1})$ with $k\geq 1$ and $A_k\not=0$, 
then equation (1.2) has two solutions,
$$
\leqalignno{
	F_1(x)&={-B+\sqrt{B^2-4AC}\over 2C}; \cr
	F_2(x)&={-B-\sqrt{B^2-4AC}\over 2C}. \cr
}
$$
Note that $F_1(x)=-A_kx^k + O(x^{k+1})$ and 
$F_2(x)=-1/C(0) + O(x)$.
\smallskip\noindent
(iii.1) In the case of $F_1(x)$, let 
$$F_1(x)={-A_kx^k \over D(x)-x^{k+2} G(x)}.$$
Then, $G(x)$ satisfies (3.7) with polynomials
$A^*, B^*, C^*$ defined by (3.5) (see the proof of Lemma 3.2). 
Since $B^*(0)=1, C^*(0)=0, C^*(x)\not=0$,
the $H$-fraction expansion of $G(x)$ exists and 
is ultimately periodic by case (i), so does the $H$-fraction expansion of $F_1(x)$.
\smallskip\noindent
(iii.2) In the case of $F_2(x)$, let 
$$F_1(x)={-1/C(0) \over D(x)-x^{2} G(x)}.\leqno{(3.12)}$$
Then, $G(x)$ satisfies (3.7) with polynomials
$A^*, B^*, C^*$ defined (same proof as Lemma 3.2):
$$
\leqalignno{
	A^*(x) &= \bigl(D^2AC(0)-BD +C/C(0)\bigr)/x^{2}; \cr
	B^*(x) &= -2ADC(0) + B ; &(3.13) \cr
	C^*(x) &= C(0)Ax^{2}. \cr
}
$$
Since $B^*(0)=1, C^*(0)=0, C^*(x)\not=0$,
the $H$-fraction expansion of $G(x)$ exists and 
is ultimately periodic by case (i), so does the $H$-fraction expansion of $F_2(x)$.
\smallskip\noindent
(iv) If $B(x)=0$,  $C(0)=1$ (or $C(0)\not=0$)
and $A(x)= -(a_k x^k)^2 + O(x^{2k+1})$ for some $k\in \setN$ and $a_k\not=0$,
then $F(x)$ exists
$$
F(x)=\sqrt{-A(x)\over C(x)}
=\sqrt{(a_k x^k)^2 + \cdots\over C(x)}
=a_kx^k\sqrt{1 + \cdots \over C(x)}
$$
Let 
$$F(x)={a_kx^k \over D(x)-x^{k+2} G(x)}.$$
Then, $G(x)$ satisfies (3.7) with 
$A^*, B^*, C^*$ defined (same proof as Lemma 3.2):
$$
\leqalignno{
	A^*(x) &= (D^2A +C a_k^2x^{2k}) / x^{3k+2}; \cr
	B^*(x) &= -2ADx^{k+2} / x^{3k+2}  ; &(3.14) \cr
	C^*(x) &= Ax^{2k+4} / x^{3k+2}. \cr
}
$$
If $p\not=2$, then  $A^*, B^*,C^*$ are polynomials
such that $B^*(0)\not=0, C^*(0)=0, C^*(x)\not=0$.
The $H$-fraction expansion of $G(x)$ exists and 
is ultimately  periodic by case (i), so does the $H$-fraction expansion of $F_2(x)$. 
\medskip
The periodicity of the Hankel determinant sequece $H(F)$ is a consequence 
of Lamma 3.4 stated below.\qed
\medskip

\proclaim Lemma 3.4.
If the $H$-fraction expansion of a power series $F$ is
ultimately periodic, then the Hankel determinant sequece
$H(F)$ is ultimately periodic.

{\it Proof}.
Using the notations of Theorem 2.1 the two sequences $(v_i)$ and $(k_i)$
can be written as 
$$
\leqalignno{
	(v_i) &= (v_0, v_1, \ldots, v_{m-1}, (v_m, v_{m+1}, \ldots, v_{m+t-1})^*);\cr
	(k_i) &= (k_0, k_1, \ldots, k_{m-1}, (k_m, k_{m+1}, \ldots, k_{m+t-1})^*).\cr
}
$$
Let 
$$
\leqalignno{
	\gamma_1&=\prod_{i=m}^{m+t-1} (-1)^{k_i(k_i+1)/2}, \quad 
	\gamma_2=\prod_{i=m}^{m+t-1} v_i^{s_m-s_i}, \quad
	\gamma_3=\prod_{i=0}^{m-1} v_i,  \cr
	\beta&=\prod_{i=m}^{m+t-1} v_i, \quad
	\gamma=\gamma_1 \gamma_3^r \gamma_2 \beta^{r-s_m}, \cr
	r&=s_{m+t} - s_m, \quad
	\eta= \lceil{\ell-m\over t} \rceil, \quad \rho=\ell -m - \eta t.
}
$$
For each $\ell\geq m$ we have
$$
H_{s_\ell}(F)= \prod_{i=0}^{\ell-1} (-1)^{k_i(k_i+1)/2} v_i^{s_\ell-s_i}
$$
and
$$
H_{s_{\ell}+r}(F)= \prod_{i=0}^{\ell+t-1} (-1)^{k_i(k_i+1)/2} v_i^{s_\ell+r-s_i}.\leqno{(3.15)}
$$
If $p=2$, then $H_{s_\ell}(F)=H_{s_\ell+r}(F)=1$. Hence $H(F)$ is 
ultimately periodic. For general $p$ we need evaluate (3.15).
$$
\leqalignno{
	H_{s_{\ell}+r}(F)
	&= 
	\prod_{i=0}^{\ell-1} (-1)^{k_i(k_i+1)/2} v_i^{s_\ell+r-s_i}\times
	\prod_{i=\ell}^{\ell+t-1} (-1)^{k_i(k_i+1)/2} v_i^{s_\ell+r-s_i}\cr
	&=  H_{s_\ell}(F)
	\prod_{i=0}^{\ell-1} v_i^{r}\times \gamma_1
	\prod_{i=\ell}^{\ell+t-\rho-1} v_i^{s_\ell+r-s_i}
	\prod_{i=\ell+t-\rho}^{\ell+t-1} v_i^{s_\ell+r-s_i} \cr
&=  H_{s_\ell}(F)
	\prod_{i=0}^{\ell-1} v_i^{r}\times \gamma_1
	\prod_{i=\ell}^{\ell+t-\rho-1} v_i^{r}
	\prod_{i=\ell}^{\ell+t-\rho-1} v_i^{s_\ell-s_i}
	\prod_{i=\ell-\rho}^{\ell-1} v_i^{s_\ell-s_i} \cr
&=  H_{s_\ell}(F) \gamma_1
	\prod_{i=0}^{\ell+t-\rho-1} v_i^{r}
	\prod_{i=\ell}^{\ell+t-\rho-1} v_i^{s_\ell-s_i}
	\prod_{i=\ell-\rho}^{\ell-1} v_i^{s_\ell-s_i} \cr
&=  H_{s_\ell}(F) \gamma_1 \gamma_3^r
	\prod_{i=m}^{\ell+t-\rho-1} v_i^{r}
\prod_{i=\ell-\rho}^{\ell+t-\rho-1} v_i^{s_\ell-s_i} \cr
&=  H_{s_\ell}(F) \gamma_1 \gamma_3^r
	\prod_{i=m}^{m+\eta t+t-1} v_i^{r}
	\prod_{i=m}^{m+t-1} v_i^{s_{m+\rho}-s_i} \cr
&=  H_{s_\ell}(F) \gamma_1 \gamma_3^r
	\prod_{i=m}^{m+\eta t+t-1} v_i^{r}
	\prod_{i=m}^{m+t-1} v_i^{s_{m+\rho}-s_m} \times \gamma_2\cr
&=  H_{s_\ell}(F) \gamma_1 \gamma_3^r \gamma_2
	\beta^{(\eta +1)r}
\beta^{s_{m+\rho} -s_m}
	\cr
&=  H_{s_\ell}(F) 
	\beta^{s_\ell} \gamma. &(3.16)\cr
}
$$
We apply (3.16) recursively and get
$$
H_{s_\ell+2\pi r}(F) 
= H_{s_\ell}(F) \beta^{2\pi s_\ell} \beta^{\pi(2\pi-1)} \gamma^{2\pi}.
$$
Choose $\pi$ such that $\beta^\pi=1$ and $\gamma^{2\pi}=1$.  Then
$$
H_{s_\ell+2\pi r}(F) 
= H_{s_\ell}(F).
$$
So that $H(F)$ is ultimately periodic.\qed

\medskip

The following notations are used for continued fraction 
$$
\KK_{n=0}^\infty {v_n\over u_{n+1}}
=
{v_0\over u_1} 
{\atop +} { v_1\over u_2}
{\atop +} { v_2\over u_3}
{\atop +} { v_3\over u_4} 
{\atop +} {\cdots} 
={v_0\hfill\over {
		u_1+\displaystyle{\strut v_1\hfill\over{
		u_2+\displaystyle{\strut v_2\hfill\over{
		u_3+\displaystyle{\strut v_3\hfill\over{
		u_4+\displaystyle{ \ddots
}}}}}}}}}.
$$

{\it Example (i.1).} Let $p=5$ and 
$$
F={1-\sqrt{1-{4x\over 1-x^4}} \over 2x}\in\setF_5[[x]]
$$
or
$$
-1 +  (1-x^4)F + (-x+x^5)F^2=0.
$$
By Algorithm 3.3 [HFrac], the power series $F$ has the following
$H$-fraction expansion
$$
\leqalignno{
F= {1 \over 1+4x }
\Kadd+
\Bigl(&{4x^2 \over 1+3x} 
\Kadd+{ 3x^2 \over 1+x }
\Kadd+{ 4x^3 \over 1+3x+2x^2}  \cr
&\quad 
\Kadd+{ 4x^3 \over 1+x} 
\Kadd+{ 3x^2 \over 1+3x }
\Kadd+{ 4x^2 \over 1+3x} 
\Kadd+{4x^2 \over 1+3x}\Kadd+ \Bigr)^*.\cr
}
$$
Hence, $m=1, t=7$, 
$(k_i)_{i\geq 0}=   (0, (0, 0, 1, 0, 0, 0, 0)^*)$,
$$(s_i)_{i\geq 0}= (0,1,  2, 3, 5, 6, 7, 8, 9, 10,\ldots),$$
$r=s_{m+t}-s_m=8$, 
and $\beta=4, \gamma_1=-1, \gamma_2=4, \gamma_3=1, \gamma=4, \pi=2$.
So that the period is less than or equal to $2\pi r=32$, starting before 
the index $m=1$.
Checking the first $32+1=33$ terms,  the period is equal to 16,
starting from index 0.
Hence
$$
H(g)= (1, 1, 1, 2, 0, 2, 4, 1, 4, 1, 4, 2, 0, 2, 1, 1)^* .
$$

{\it Example (i.2).} Let $p=2$ and 
$$
F={1-\sqrt{1-{4x\over 1-x^4}} \over 2x} \in\setF_2[[x]] \leqno{(3.17)}
$$
or
$$
-1 +  (1-x^4)F + (-x+x^5)F^2=0. \leqno{(3.18)}
$$
By Algorithm 3.3 [HFrac], we get the following $H$-fraction  expansion
$$
\leqalignno{
F= {1 \over 1+x }
\Kadd+
\Bigl(&{x^2 \over 1} 
\Kadd+{ x^4 \over 1 }
\Kadd+{ x^6 \over 1}  
\Kadd+{ x^4 \over 1 }
\Kadd+{ x^2 \over 1} 
\Kadd+{x^2 \over 1}\Kadd+ \Bigr)^*.\cr
}
$$
Hence, $m=1, t=6$, 
$(k_i)_{i\geq 0}= (0, (0, 2,  2, 0, 0, 0)^*)$,
$$(s_i)_{i\geq 0}= (0, 1, 2,  5, 8, 9, 10,11,\ldots),$$
$r=s_7-s_1=10$.
The period is less than or equal to 10, starting before the index $s_m=1$.
Checking the first $10+1=11$ terms in 
$H(f)$, which are
$(1,1, 1, 0, 0, 1, 0, 0, 1, 1, 1,   \ldots),
$
we see that the period is equal to $10$ and
$H(f)=
(1,1, 1, 0, 0, 1, 0, 0, 1, 1)^*.
$
\medskip

{\it Example (iii.1).} Let $p=2$ and  $G=xF$ where $F$ is defined in Example (i.2) by (3.17) or (3.18). We have
$$
G={1-\sqrt{1-{4x\over 1-x^4}} \over 2} \in\setF_2[[x]]
$$
and
$$
-x +  (1-x^4)G + (-1+x^4)G^2=0 \hbox{\quad \rm with\quad} G(0)=0. \leqno{(3.19)}
$$
Since the coefficient of $G^2$ has constant term, we cannot apply Algorithm 3.3 directly.
Let
$$
G={x\over 1+x+x^2+x^3 G_1}.
$$
Equation (3.19) becomes
$$
x^3 +  (1+x^4)G_1 + x^3G_1^2=0. 
$$
By Algorithm 3.3 [HFrac], we get the following $H$-fraction  expansion
$$
\leqalignno{
G_1= {x^3 \over 1+x^4 }
\Kadd+
\Bigl(&{x^6 \over 1} 
\Kadd+{ x^4 \over 1+x^2 }
\Kadd+{ x^4 \over 1 }
\Kadd+{ x^6 \over 1+x^4} 
\Kadd+ \Bigr)^*.\cr
}
$$
Hence
$$
G= 
{x\over 1+x+x^2}
\Kadd+
\Bigl(
{x^6 \over 1+x^4 }
\Kadd+ {x^6 \over 1} 
\Kadd+{ x^4 \over 1+x^2 }
\Kadd+{ x^4 \over 1 }
\Kadd+ \Bigr)^*.
$$

{\it Example (iii.2).} Let $p=2$ and  $$
G={1+\sqrt{1-{4x\over 1-x^4}} \over 2} \in\setF_2[[x]]
$$
and
$$
-x +  (1-x^4)G + (-1+x^4)G^2=0 \hbox{\quad \rm with\quad} G(0)=1. \leqno{(3.20)}
$$
Since the coefficient of $G^2$ has constant term, we cannot apply Algorithm 3.3 directly.
Let
$$
G={1\over 1+x+x^2 G_1}.
$$
Equation (3.20) becomes
$$
(x+x^3) +  (1+x^4)G_1 + x^3G_1^2=0. 
$$
By Algorithm 3.3 [HFrac], we get the following $H$-fraction  expansion
$$
G_1= {x \over 1+x^2 }
\Kadd+
\Bigl({x^4 \over 1} 
\Kadd+{ x^6 \over 1+x^4 }
\Kadd+{ x^6 \over 1 }
\Kadd+{ x^4 \over 1+x^2} 
\Kadd+ \Bigr)^*.
$$
Hence
$$
G= 
{1\over 1+x}
\Kadd+
{x^3 \over 1+x^2 }
\Kadd+
\Bigl({x^4 \over 1} 
\Kadd+{ x^6 \over 1+x^4 }
\Kadd+{ x^6 \over 1 }
\Kadd+{ x^4 \over 1+x^2} 
\Kadd+ \Bigr)^*.
$$

{\it Example (iv).} Let $p=3$ and  
$$
F=\sqrt{x^2-x^3\over 1+x^3}\in\setF_3[[x]]
$$
or
$$
(-x^2+x^3)+(1+x^3)F^2=0
$$
Since the coefficient of $F$ is zero, we cannot apply Algorithm 3.3 directly.
Let
$$
F={x\over 1+2x+ x^3 F_1}.
$$
Equation (3.20) becomes
$$
2 +  (1+x+x^2)F_1 + (2x^3+x^4)F_1^2=0. 
$$
By Algorithm 3.3 [HFrac], we get the following $H$-fraction  expansion
$$
F_1= {1 \over 1+x }
\Kadd+
\Bigl({x^2 \over 1+x} 
\Kadd+{ 2x^2 \over 1+x }
\Kadd+{ x^2 \over 1+x }
\Kadd+{ 2x^3 \over 1+2x }
\Kadd+{ 2x^3 \over 1+x} 
\Kadd+ \Bigr)^*.
$$
Hence, the $H$-fraction expansion of $F$ is 
$$
{x \over 1+2x }
\Kadd+
{x^3 \over 1+x }
\Kadd+
\Bigl({x^2 \over 1+x} 
\Kadd+{ 2x^2 \over 1+x }
\Kadd+{ x^2 \over 1+x }
\Kadd+{ 2x^3 \over 1+2x }
\Kadd+{ 2x^3 \over 1+x} 
\Kadd+ \Bigr)^*.
$$

\medskip

The proof of Theorem 1.1 is also valid for super $1$-fraction.

\proclaim Theorem 3.5.
Let $p$ be a prime number and $F(x)\in \setF_p[[x]]$ be a power series
satisfying the following quadratic functional equation
$$A(x)+B(x)F(x)+ C(x)F(x)^2 = 0, $$
where $A(x), B(x), C(x)\in \setF_p[x] $ are three polynomials with one of the following conditions 
\smallskip
(i) $B(0)=1, \ C(0)=0,\ C(x)\not=0$;
\smallskip
(ii) $B(0)=1,\ C(x)=0$;
\smallskip
(iii) $B(0)=1,\ C(0)\not=0,\ A(0)=0$;
\smallskip
(iv) $B(x)=0,\  C(0)=1, \ A(x)= -(a_k x^k)^2 + O(x^{2k+1})$ for some $k\in \setN$ and $a_k\not=0$ when $p\not=2$.
\smallskip\noindent
Then, the super $1$-fraction expansion of $F(x)$ exists and 
is ultimately periodic. 


\section{4. Application to automatic sequences } 
\def\sec{4}

{\it Proof of Theorem 1.2}. 
Let $f(x)=G_{a,b}(x)\in\setF_2[[x]]$. Then
$$
\leqalignno{
	x^{2^a}f(x)&= \sum_{n=0}^\infty {x^{2^{n+a}} \over 1-x^{2^{n+b}}};\cr
	x^{2^{a+1}}f(x^2)&= \sum_{n=1}^\infty {x^{2^{n+a}} \over 1-x^{2^{n+b}}};\cr
	x^{2^{a}}f(x^2)&=f(x)-{1\over 1-x^{2^b}}; \cr
	1 +  (1+x^{2^b})f(x)& + x(1+x^{2^b})x^{2^a-1}f(x)^2= 0. \cr
}
$$
The above equation is of type (1.2). By Theorem 1.1 the Hankel 
determinant sequence $H(f)$ is ultimately periodic.\qed
\medskip

The following corollary is obtained by Algorithm 3.3. 
The case for
the regular paperfolding sequence, i.e., $a=0, b=2$,
is verified in
[Section 3, Example (i.2)].

\proclaim Corollary 4.1.
Let $G_{a,b}(x)$ be power series in $\setF_2[[x]]$ defined by (1.3).
Over the field $\setF_2$ we have
$$
\leqalignno{
	H(G_{0,0})&= (1)^* ; \cr
	H(G_{0,1})&= 1,1,(0)^* ; \cr
	H(G_{1,0})&= (1)^* ; \cr
	H(G_{0,2})&= (1,1,1,0,0,1,0,0,1,1)^* ; \cr
	H(G_{1,1})&= (1,1, 0, 0, 1, 1)^* ; \cr
H(G_{2,0})&= (1,1, 0, 0)^* ; \cr
	H(G_{0,3})&= ({1}^{5}{0}^{2}{1}^{1}{0}^{6}{1}^{3}{0}^{2}{1}^{2}{0}^{2}{1}^{2}{0}^{4}{1}^{1}{0}^{4}{1}^{1}{0}^{2}{1}^{1}{0}^{2}{1}^{1}\cr
	& \qquad {0}^{4}{1}^{1}{0}^{4}{1}^{2}{0}^{2}{1}^{2}{0}^{2}{1}^{3}{0}^{6}{1}^{1}{0}^{2}{1}^{4})^* ; \qquad \hbox{\rm [period is 74]} \cr
	\noalign{\goodbreak}
	H(G_{1,2})&= 1,1,1,(0)^* ;  \cr
	H(G_{2,1})&= (1,1, 1, 1, 1, 1, 0, 0)^* ; \cr
	H(G_{3,0})&= (1,1, 0, 0, 0, 0, 0, 0)^* ; \cr
	H(G_{0,4})&= (1^90^21^10^2\cdots 1^10^21^8)^*  \quad \hbox{\rm [period is 1078]} ; \cr
	H(G_{1,3})&= (1,1, 1, 1, 1, 0, 0, 0, 0, 1, 1, 0, 0, 0, 0, 1, 1, 1, 1, 1)^* ; \cr
	H(G_{2,2})&= (1,1, 0, 0, 0, 0, 0, 0, 1, 1, 0, 0)^* ; \cr
	H(G_{3,1})&= (1,1, 1, 0, 0, 0, 0, 1, 1, 1, 0, 0, 0, 0, 0, 0)^* ; \cr
	H(G_{4,0})&= (1,1, 0, 0, 0, 0, 0, 0, 0, 0, 0, 0, 0, 0, 0, 0)^* . \cr
}
$$

\medskip


{\it Proof of Proposition 1.3}.
Let $f(x)=\sum_{n\geq 0} u_n x^n\in\setF_2[[x]]$ where $(u_n)$ is the Rudin-Shapiro sequence 
defined by (1.4). Then 
$$
x^3+(1+x)^4f(x) + (1+x)^5f^2(x) =0.  \leqno{(4.1)}
$$
Since $u_0=u_1=u_2=0$, 
$$
	f_1(x)=f(x)/x  ; \quad
	f_2(x)=f(x)/x^2; \quad
	f_3(x)=f(x)/x^3.
$$
From (4.1) we derive
$$
\leqalignno{
	x^2+(1+x)^4f_1(x) + (1+x)^5xf_1^2(x) &= 0 ;  \cr 
	x+(1+x)^4f_2(x) + (1+x)^5x^2f_2^2(x) &= 0 ;  \cr 
	1+(1+x)^4f_3(x) + (1+x)^5x^3f_3^2(x) &= 0 .  \cr 
}
$$
Theorem 1.3 follows from Algorithm 3.3.\qed
\medskip

{\it Proof of Proposition 1.4}.
It is well known that [BV13]
$$
S(x)=(1+x+x^2)S(x^2) \in\setQ[[x]]. \leqno{(4.2)}
$$
Since $S(x) \pmod 2$  is rational, there exists a positive integer $N$ such that $H_k(S)\equiv 0\pmod 2$ for all $k\geq N$. We must use the grafting technique,
introduced in [Ha13, Section 2].
First, the $H$-fraction of $S(x)$ is
$$
S(x)={1\over {1-x-
\displaystyle{x^2\over 1 +2x +
\displaystyle{2x^2 \over 1 - 
\displaystyle{2x^3\over 1-3/2 x + 11/4 x^2 + \ddots}}  }}}.
$$
The even number 2 occurs in the 
sequence $(v_j)$, in particular at
position $v_2$ in view of (2.1). Define $G(x)$ by 
$$
S(x) =
{1\over {1-x-
\displaystyle{x^2\over 1 +2x+ 
	\displaystyle{2x^2 G(x) 
}  }}}.
$$
From (4.2) the power series $G(x)$ satisfies 
the following relation
$$
(1+x+x^2) + (1+x+x^2) G(x) + x^4 G(x^2) \equiv 0 \pmod 2.
$$
By Algorithm 3.3 we get $H(G)\equiv (1, 1, 0, 0)^* \pmod 2$. 
By Lemma~3.2, $H_n(S)=(-1)^n2^{n-2} H_{n-2}(G)$.
Hence 
$$H_n(S)/2^{n-2} \equiv (0,0,1,1)^* \pmod 2.$$
\medskip
In the same manner, $B(x)$ is a rational function modulo $2$. We use the grafting technique. Since [BV13]
$$
B(x)=2-(1+x+x^2) B(x^2) \leqno{(4.3)}
$$
and
$$
B(x)={1\over
	1+x+\displaystyle{x^2\over{
			1+\displaystyle{2 x^2 \over {
					1-2x+\displaystyle{x^4\over \ddots}
				}
}}}},
$$
we define $U(x)$ by
$$
B(x)={1\over
	1+x+\displaystyle{x^2\over{
1+{2 x^2 U(x)}}}}.
$$
From (4.3) the power series $U(x)$ satisfies the following fonctional equation 
$$
(1+x+x^2) + (1+x+x^2)U(x) + x^4 U(x^2) \equiv 0 \pmod 2.
$$
By Algorithm 3.3 we get 
$H(U)=(1,1,0,0)^* \pmod 2$.
On the other hand,
the Hankel determinant
$
H_n(B) = -2^{n-2} H_{n-2}(U)
$
by Lemma 3.2.
Hence,
$H_n(B)/2^{n-2} \equiv (0,0,1,1)^* \pmod 2.$\qed

\bigskip


\goodbreak
\bigskip
{\bf Acknowledgements}. The author should like to thank Zhi-Ying Wen who 
suggested that I study
the Hankel determinants of the Thue-Morse sequence back to 1991,
and who invited me to Tsinghua University where the paper is finalized.
The author also thanks 
Yann Bugeaud and Jia-Yan Yao for valuable discussions.  

\bigskip

\vskip 2cm
\bigskip
\centerline{\bf References}
\bigskip

{
\eightpoint

\article Al87|Allouche, J.-P|Automates finis en th\'eorie des nombres|Expo. Math.|5|1987|239--266|

\divers Al12|Allouche, J.-P|On the Stern sequence and its twisted version, {\sl Integers}, {\bf 12} ({\oldstyle 2012}), A58|

\article APWW|Allouche, J.-P.; Peyri\`ere, J.; Wen, Z.-X.; Wen, Z.-Y|
Hankel determinants of the Thue-Morse sequence|Ann. Inst. Fourier, 
Grenoble|48|1998|1--27|

\divers Ba10|Bacher, Roland|Twisting the Stern sequence, {\tt arxiv.org/abs/1005.5627}, {\oldstyle 2010}, 19 pages|

\divers BH13|Bugeaud, Yann; Han, Guo-Niu|A combinatorial proof of the non-vanishing of Hankel determinants of the Thue--Morse sequence, 16 pages, 2013|

\divers BHWY|Bugeaud, Yann; Han, Guo-Niu; Wen, Zhi-Xing; Yao, Jia-Yan|Hankel determinant calculus 
for automatic sequences, III: irrationality exponent, {\sl in preparation}, {\oldstyle 2013}|

\article Bu10|Buslaev, V.I|On Hankel determinnats of functions given by
their expansions in $P$-fractions|Ukrainian Math. J.|62|2010|358--372|

\article Bu11|Bugeaud, Yann|On the rational approximation to the Thue-Morse-Mahler numbers|Ann. Inst. Fourier, Grenoble|61|2011|2065--2076|

\article BV13|Bundschuh, Peter; Vaananen, Keijo|Algebraic independence of the
generating functions of Stern's sequence and of its twist|J. Th\'eorie des Nombres de Bordeaux|25|2013|43--57|

\divers Ci13|Cigler, Johann|A special class of Hankel determinants, preprint, {\oldstyle 2013}|

\article CV12|Coons, Michael; Vrbik, Paul|An irrationality measure for regular paperfolding numbers|Journal of Integer Sequences|15|2012|Article 12.1.6|

\article Co13|Coons, Michael|On the rational approximation of the sum of the reciprocals of the Fermat numbers|The Ramanujan Journal|30|2013|39--65|

\article Fl80|Flajolet, Philippe|Combinatorial aspects of continued fractions|Discrete Math.|32|1980|125--161|

\divers GWW|Guo, Yingjun; Wu, Wen; Wen, Zhixiong|On the irrationality exponent of the regular paperfolding numbers, arxiv.org/abs/1310.2138, {\oldstyle 2013}|

\divers Ha13|Han, Guo-Niu|Hankel Determinant Calculus for the Thue-Morse and related sequences, {\sl preprint}, 21 pages, {\oldstyle 2013}|

\divers Kr98|Krattenthaler, Christian|Advanced determinant calculus,
{\sl S\'em. Lothar. Combin.},
{\bf B42q} ({\oldstyle1998}), 67pp|

\article Kr05|Krattenthaler, Christian|Advanced determinant calculus: A complement|Linear Algebra and its Applications|411|2005|68--166|

\article St58|Stern, M.A.|\"Uber eine zahlentheoretische Funktion|J. Reine Angew. Math|55|1858|193--220|

\divers Vi83|Viennot, X|Une th\'eorie combinatoire des polyn\^omes 
orthogonaux g\'en\'eraux, {\sl UQAM, Montreal, Quebec}, {\oldstyle 1983}|

\livre Wa48|Wall, H. S|Analytic theory of continued fractions|Chelsea publishing company, Bronx, N.Y., {\oldstyle 1948}|

\divers WiRP|Wikipedia|Regular paperfolding sequence, {\it revision September 30, 2013}|

\bigskip

}
\bigskip\bigskip
\hbox{
\vtop{\halign{#\hfil\cr
I.R.M.A. UMR 7501\cr
Universit\'e de Strasbourg et CNRS\cr
7, rue Ren\'e-Descartes\cr
F-67084 Strasbourg, France\cr
\noalign{\smallskip}
{\tt guoniu.han@unistra.fr}\cr}}}


\end